%


\magnification=\magstep1
\input amstex
\documentstyle{amsppt}
\NoRunningHeads

\def\pf{{\hfill$\square$}}
\def\c{\cite}
\def\st{such that }
\def\iff{if and only if }
\def\acs{atomless complete subalgebra }
\def\calp{\Cal P}
\def\ind{\text{ind}}
\def\os{\overline{s}}
\def\da{\dot{a}}
\def\dx{\dot{x}}
\def\dox{\dot{X}}
\topmatter

\title
A complete Boolean algebra that has no proper atomless complete sublagebra
\endtitle

\author
Thomas Jech and Saharon Shelah
\endauthor

\affil
The Pennsylvania State University,  \\ 
The Hebrew University and Rutgers University
\endaffil 

\thanks
The first author was supported in part by an NSF grant
DMS-9401275.
\endgraf
The second author was partially supported by the U.S.--Israel
Binational Science Foundation. Publication No. 566 \endthanks
\address
Department of Mathematics, The Pennsylvania State University,
University Park, PA 16802, USA
\endgraf
School of Mathematics, The Hebrew University, Jerusalem, Israel, and
Department of Mathematics, Rutgers University, New Brunswick, NJ 08903, USA
\endaddress
\email jech\@math.psu.edu, shelah\@math.huji.ac.il \endemail

\abstract
There exists a complete atomless Boolean algebra that has no proper
\acs\!\!.
\endabstract

\endtopmatter

\document

\baselineskip 18pt

An atomless complete Boolean algebra $B$ is {\it simple} \c5 if it has
no \acs $A$ \st $A\neq B$.  The question whether such an algebra $B$
exists was first raised in \c8 where it was proved that $B$ has no
proper \acs \iff $B$ is {\it rigid} and {\it minimal}.  For more on
this problem, see \c4, \c5 and \c{1, p. 664}.

Properties of complete Boolean algebras correspond to properties of
generic models obtained by forcing with these algebras.  (See \c6, 
pp. 266--270; we also follow \c6 for notation and terminology of
forcing and generic models.)  When in \c7 McAloon constructed a generic
model with all sets ordinally definable he noted that the corresponding
complete Boolean algebra is {\it rigid}, i.e. admitting no nontrivial
automorphisms.  In \c9 Sacks gave a forcing construction of a real number
of minimal degree of constructibility.  A complete Boolean algebra $B$
that adjoins a minimal set (over the ground model) is {\it minimal} in 
the following sense:
$$
\split
\text{If $A$ is a complete atomless subalgebra of $B$ then there 
exists}\\ \text{a partition $W$ of $1$ \st for every $w\in W$, 
$A_w=B_w$,}\\ \text{where $A_w=\{a\cdot w : a\in A\}.$}\endsplit\tag1
$$ 

In \c{3}, Jensen constructed, by forcing over $L$, a definable real
number of minimal degree.  Jensen's construction thus proves that in
$L$ there exists rigid minimal complete Boolean algebra.  This has
been noted in \c8 and observed that $B$ is rigid and minimal \iff it
has no proper \acs\!\!.  McAloon then asked whether such an algebra
can be constructed without the assumption that $V=L$.  In \c5 simple
complete algebras are studied systematically, giving examples (in $L$)
for all possible cardinalities.

In \c{10} Shelah introduced the $(f,g)$-bounding property of forcing
and in \c2 developed a method that modifies Sacks' perfect tree forcing
so that while one adjoins a minimal real, there remains enough freedom
to control the $(f,g)$-bounding property.  It is this method we use
below to prove the following Theorem:

\proclaim
{Theorem}  There is a forcing notion $\calp$ that adjoins a real number
$g$ minimal over $V$ and \st $B(\calp)$ is rigid.
\endproclaim

\proclaim
{Corollary}  There exists a countably generated simple complete 
Boolean algebra.
\endproclaim

The forcing notion $\calp$ consists of finitely branching perfect 
trees of height $\omega$.  In order to control the growth of trees $T\in\calp$,
we introduce a {\it master tree} $\Cal T$ \st every $T\in \calp$ will be a
subtree of $\Cal T$.  To define $\Cal T$, we use the following fast growing
sequences of integers $(P_k)_{k=0}^{\infty}$ and $(N_k)_{k=0}^{\infty}$:
$$
P_0=N_0=1,\quad  P_k=N_0\cdot\dots\cdot N_{k-1}, \quad N_k=2^{P_k}\tag2
$$
(Hence $N_k=1,2,4,256,2^{2^{11}},\dots$).
\proclaim
{Definition} {\rm The {\it master tree} $\Cal T$ and the {\it index function}
ind:
\roster
\item"{(3)(i)}" $\Cal T\subset [\omega]^{<\omega}$,
\item"{(ii)}" ind is a one-to-one function of $\Cal T$ onto $\omega$,
\item"{(iii)}" ind $(< \ >)=0$,
\item"{(iv)}" if $s,t\in \Cal T$ and length(s) $<$ length(t) then
$\ind(s)< \ind(t),$
\item"{(v)}" if $s,t\in \Cal T$, length(s) $=$ length(t) and $s<_{lex}t$ then
$\ind(s)< \ind(t),$
\item"{(vi)}" if $s\in \Cal T$ and $\ind(s)=k$ then $s$ has exactly $N_k$
successors in $\Cal T$, namely all $s^{\frown}i$, $i=0,\dots,N_k-1$.
\endroster
\endproclaim

The forcing notion $\calp$ is defined as follows:
\proclaim
{Definition} {\rm $\calp$ is the set of all subtrees $T$ of $\Cal T$ that
satisfy the following:
$$
\split
\text{for every $s\in T$ and every $m$ there exists some $t\in T$,
$t\supset s$,}\\ \text{ \st $t$ has at least $P_{\ind(t)}\!^m$ successors in $T$.}
\endsplit\tag4 
$$
\endproclaim

(We remark that $\Cal T\in \calp$ because for every $m$ there is a $K$ \st
for all $k\geq K$, $P_k\!^m\leq 2^{P_k}=N_k$.)  

When we need to verify
that some $T$ is in $\calp$ we find it convenient to replace (4) by
an equivalent property:
\proclaim
{Lemma}  A tree $T\subseteq J$ satisfies {\rm (4)} \iff
\roster
\item"{(5)(i)}" every $s\in T$ has at least one successor in $T$,
\item"{(ii)}"  for every $n$, if $\ind(s)=n$ and $s\in T$ then there exists
a $k$ \st if $\ind(t)=k$ then $t\in T$, $t\supset s$ and $t$ has at least
$P_k\!^n$ successors in $T$.
\endroster
\endproclaim

\demo
{Proof}  To see that (5) is sufficient, let $s\in T$ and let $m$ be
arbitrary.  Find some $\os\in T$ \st $\os\supset s$ and $\ind(\os)\geq m$,
and apply (5ii).   \pf
\enddemo

The forcing notion $\calp$ is partially ordered by inclusion.  A standard 
forcing argument shows that if $G$ is a generic subset of $\calp$ then
$V[G]=V[g]$ where $g$ is the {\it generic branch}, i.e. the unique
function $g:\omega\to \omega$ whose initial segments belong to all $T\in G$.  
We shall prove that the generic branch is minimal over $V$, and that the
complete Boolean algebra $B(\calp)$ admits no nontrivial automorphisms.

First we introduce some notation needed in the proof:
$$
\text{For every $k$, $s_k$ is the unique $s\in \Cal T$ \st $\ind(s)=k$.}\tag6
$$
$$
\text{If $T$ is a tree then $s\in\text{trunk}(T)$ if for all $t\in T$,
either $s\subseteq t$ or $t\subseteq s$.}\tag7
$$
$$
\text{If $T$ is a tree and $a\in T$ then 
$(T)_a=\{s\in T : s\subseteq a\text{ or } a\subseteq s\}. $}\tag8
$$

Note that if $T\in\calp$ and $a\in T$ then $(T)_a\in\calp$.  We shall
use repeatedly the following technique:

\proclaim
{Lemma} Let $T\in\calp$ and, let $l$ be an integer and let $U=T\cap
\omega^{l}$ (the $l^{\text{th}}$ level of $T$).  Let $\dx$ be a name
for some set in $V$.  For each $a\in U$ let $T_a\subseteq(T)_a$ and
$x_a$ be \st $T_a\in \calp$ and $T_a\Vdash\dx=x_a$.

Then $T'=\bigcup\{T_a:a\in U\}$ is in $\calp$, $T'\subseteq T$,
$T'\cap \omega^{l}=T\cap \omega^{l}=U$, and $T'\Vdash \dx\in\{x_a:a\in U\}$.
\pf
\endproclaim

We shall combine this with {\it fusion}, in the form stated below:

\proclaim
{Lemma}  Let $(T_n)_{n=0}^{\infty}$ and $(l_n)_{n=0}^{\infty}$ be
\st each $T_n$ is in $\calp$, $T_0\supseteq T_1\supseteq\dots\supseteq
T_n\supseteq\dots$, $l_0<l_1<\dots<l_n<\dots$, $T_{n+1}\cap
\omega^{l_n}=$ $T_n\cap \omega^{l_n}$, and \st
$$
\split
\text{for every $n$, if $s_n\in T_n$ then there exists some
$t\in T_{n+1}$, $t\supset s_n$, with}\\ \text{length$(t)<l_{n+1}$, \st
$t$ has at least $P_{\ind(t)}\!^n$ successors in $T_{n+1}$.}\endsplit\tag9
$$
Then $T=\bigcap^{\infty}_{n=0}T_n\in\calp$.
\endproclaim

\demo
{Proof}  To see that $T$ satisfies (5), note that if $s_n\in T$ 
then $s_n\in T_n$, and the node $t$ found by (9) belongs to $T$.  \pf
\enddemo

We shall now prove that the generic branch is minimal over $V$:

\proclaim
{Lemma}  If $X\in V[G]$ is a set of ordinals, then either $X\in V$
or $g\in V[X]$.
\endproclaim

\demo
{Proof} The proof is very much like the proof for Sacks' forcing.  Let
$\dox$ be a name for $X$ and let $T_0\in\calp$ force that $\dox$ is
not in the ground model.  Hence for every $T\leq T_0$ there exist $T'$,
$T''\leq T$ and an ordinal $\alpha$ \st $T'\Vdash\alpha\in\dox$ and
$T''\Vdash\alpha\notin\dox$.  Consequently, for any $T_1\leq T$ and
$T_2\leq T$ there exist $T_1'\leq T_1$ and $T_2'\leq T_2$ and an
$\alpha$ \st both $T_1'$ and $T_2'$ decide ``$\alpha\in\dox$'' and
$T_1'\Vdash\alpha\in\dox$ \iff $T_2'\Vdash\alpha\notin\dox$.  

Inductively, we construct $(T_n)_{n=0}^{\infty}$,
$(l_n)_{n=0}^{\infty}$, $U_n=T_n\cap \omega^{l_n}$, and ordinals
$\alpha(a,b)$ for all $a,b\in U_n$, $a\neq b$, \st
\roster\item"{(10)(i)}" $T_n\in\calp$ and $T_0\supseteq T_1
\supseteq\dots\supseteq T_n\supseteq\dots$,
\item"{(ii)}" $l_0<l_1<\cdots<l_n<\cdots$,
\item"{(iii)}" $T_{n+1}\cap \omega^{l_n}=T_n\cap \omega^{l_n}=U_n$,
\item"{(iv)}" for every $n$, if $s_n\in T_n$ then there exists some
$t\in T_{n+1}$, $t\supset s_n$, with length$(t)<l_{n+1}$, \st
$t$ has at least $P_{\ind(t)}\!^n$ successors in $T_{n+1}$, 
\item"{(v)}" for every $n$, for all $a,b\in U_n,$ if $a\neq b$ then 
both $(T_n)_a$ 
and $(T_n)_b$ decide ``$\alpha(a,b)\in\dox$'' and $(T_n)_a\Vdash
\alpha(a,b)\in\dox$ \iff $(T_n)_b\Vdash\alpha(a,b)\in\dox$.
\endroster

When such a sequence has been constructed, we let
$T=\bigcap^{\infty}_{n=0} T_n$. As (9) is
satisfied, we have $T\in\calp$ and $T\leq T_0$.  If $G$ is a generic
\st $T\in G$ and if $X$ is the $G$-interpretation of $\dox$ then the
generic branch $g$ is in $V[X]$: for every $n$,
$g\upharpoonright l_n$ is the unique $a\in U_n$ with the property
that for every $b\in U_n$, $b\neq a$, $\alpha(a,b)\in X$ \iff
$(T)_a\Vdash\alpha(a,b)\in\dox$.

To construct $(T_n)_{n=0}^{\infty}$, $(l_n)_{n=0}^{\infty}$ and
$\alpha(a,b)$, we let $l_0=0$ (hence $U_0=\{s_0\}$) and proceed by
induction.  Having constructed $T_n$ and $l_n$, we first find
$l_{n+1}>l_n$ as follows:  If $s_n\in T_n$, we find $t\in T_n$,
$t\supset s_n$, \st $t$ has at least $P_{\ind(t)}\!^n$ successors in
$T_n$.  Let $l_{n+q}=$ length $(t)+1$.  (If $s_n\notin T_n$, let
$l_{n+1}=l_n+1$.)  Let $U_{n+1}=T_n\cap \omega^{l_{n+1}}$.

Next we consider, in succession, all pairs $\{a,b\}$ of district
elements of $U_{n+1}$, eventually constructing conditions
$T_a$, $a\in U_{n+1}$, and ordinals $\alpha(a,b)$, $a,b\in U_{n+1}$,
\st for all $a$, $T_a\leq (T_n)_a$ and if $a\neq b$ then either 
$T_a\Vdash\alpha(a,b)\in\dox$ and $T_b\Vdash\alpha(\da,b)\notin\dox$
or vice versa.  Finally, we let $T_{n+1}=\bigcup\{T_a:
a\in U_{n+1}\}$.

It follows that $(T_n)_{n=0}^{\infty}$, $(l_n)_{n=0}^{\infty}$ and
$\alpha(a,b)$ satisfy (10).  \pf
\enddemo

Let $B$ be the complete Boolean algebra $B(\calp)$.  We shall prove
that $B$ is rigid.  Toward a contradiction, assume that there exists
an automorphism $\pi$ of $B$ that is not the identity.  First, there
is some $u\in B$ \st $\pi(u)\cdot u=0$.  Let $p\in\calp$ be \st $p\leq
u$ and let $q\in\calp$ be \st $q\leq\pi(p)$.  Since $q\not\leq p$,
there is some $s\in q$ \st $s\notin p$.  Let $T_0=(q)_s$.

Note that for all $t\in T_0$, if $t\supseteq s$ then $t\notin p$.  
Let
$$
	A=\{\ind(t):t\in p\},
$$
and consider the following property $\varphi(x)$ (with parameters
in V):
$$
\split
\varphi(x)\;\leftrightarrow\; &\text{if $x$ is a
function from $A$ into $\omega$ then there exists }\\ &\text{a function $u$ on
$A$ whose values are finite sets of}\\ &\text{integers and for every $k\in A$,
$u(k)\subseteq N_k$ and $|u(k)|\leq P_k$,}\\
&\text{and $x(k)\in u(k)$.}\endsplit\tag11
$$
We will show that
$$p\Vdash\exists x\neg\varphi(x),\tag12$$
and
$$
\text{there exists a $T\leq T_0$ \st $T\Vdash\forall x\varphi(x)$.}\tag13
$$ 

This will yield a contradiction:  the Boolean
value of the sentence $\exists x\neg\varphi(x)$ is preserved by $\pi$,
and so 
$$T_0\leq q\leq\pi(p)\leq\pi(\|\exists x\neg\varphi(x)\|)
=\|\exists x\neg\varphi(x)\|,$$
contradicting (13).

In order to prove (12), consider the following (name for a) function
$\dx:A\to \omega$.   For every $k\in A$, let 
$$
	\dx(k)=\dot{g}(\text{length}\,(s_k)+1)\text{ if }
	s_k\subset\dot{g},\text{ and }\dx(k)=0\quad
	\text{otherwise}.
$$ 
Now if $p_1<p_2$ and $u\in V$ is a function on $A$ \st $u(k)\subseteq
N_k$ and $|u(k)|\leq P_k$ then there exist a $p_2<p_1$ and some
$k\in A$ \st $s_k\in p_2$ has at least $P_k\!^2$ successors, and
there exist in turn a $p_3<p_2$ and some $i\notin u(k)$ \st
$s_k^{\frown} i\in$~trunk$(p_3)$.  Clearly, $p_3\Vdash \dx(k)\notin
u(k)$.

Property (13) will follow from this lemma:

\proclaim
{Lemma}  Let $T_1\leq T_0$ and $\dx$ be \st $T_1$ forces that
$\dx$ is function from $A$ into $\omega$.  There exist sequences 
$(T_n)_{n=1}^{\infty}$, $(l_n)_{n=1}^{\infty}$, $(j_n)_{n=1}^{\infty}$,
$(U_n)_{n=1}^{\infty}$ and sets $z_a$, $a\in U_n,$ \st
\roster\item"{(14)(i)}" $T_n\in\calp$
and $T_1\supseteq T_2\supseteq\dots\supseteq T_n\supseteq\dots$,
\item"{(ii)}" $l_1<l_2<\dots<l_n<\dots$,
\item"{(iii)}" $T_{n+1}\cap \omega^{l_n}=T_n\cap \omega^{l_n}=U_n$,
\item"{(iv)}" for every $n$, if $s_n\in T_n$ then there exists
some $t\in T_{n+1}$, $t\supset s_n$, with length$(t)<l_{n+1}$,
\st $t$ has at least $P_{\ind(t)}\!^n$ successors in $T_{n+1}$, 
\item"{(v)}" $j_1<j_2<\dots<j_n<\dots$,
\item"{(vi)}" for every $a\in U_n$, $(T_n)_a\Vdash
\langle\dx(k):k\in A\cap j_n\rangle=z_a$,
\item"{(vii)}" for every $k\in A$, if $k\geq j_n$ then $|U_n|<P_k$,
\item"{(viii)}" for every $k\in A$, if $k<j_n$ then $|\{z_a(k):a\in
U_n\}|\leq P_k$.
\endroster
\endproclaim

Granted this lemma, (13) will follow: If we let
$T=\bigcap^{\infty}_{n=1} T_n$, then $T\in\calp$
and $T\leq T_1$ and for every $k\in A$, $T\Vdash\dx(k)\in u(k)$ where
$u(k)=\{z_a(k):a\in U_n\}$ (for any and all $n>k$).

\demo
{Proof of Lemma} We let $l_1=j_1=$ length$(s)$, $U_1=\{s\}$ and
strengthen $T_1$ if necessary so that $T_1$ decides $\langle
\dx(k):k\in A\cap j_1\rangle$, and let $z_s$ be the decided value.  We
also assume that length$(s) \geq 2$ so that $|U_1|=1<P_k$ for every
$k\in A$, $k\geq j_1$.  Then we proceed by induction.
\enddemo

Having constructed $T_n$, $l_n$, $j_n$ etc., we first find
$l_{n+1}>l_n$ and $j_{n+1}>j_n$ as follows:  If $s_n\notin
T_n$ (Case I), we let $l_{n+1}=l_n+1$ and $j_{n+1}=j_n+1$.
Thus assume that $s_n\in T_n$ (Case II).

Since length$(s_n)\leq n\leq l_n$, we choose some $v_n\in U_n$ 
\st $s_n\subseteq v_n$.  By (4) there exists some $t\in T_n$,
$t\supset v_n$, so that if ind$(t)=m$ then $t$ has at least
$P_m\!^{n+1}$ successors in $T_n$.  Moreover we choose $t$ so that
$m=\ind(t)$ is big enough so that there is at least one $k\in A$ \st
$j_n\leq k<m$.  We let $l_{n+1}=$ length$(t)+1$ and
$j_{n+1}=m=\ind(t)$.  

Next we construct $U_{n+1}, \{z_a:a\in
U_{n+1}\}$ and $T_{n+1}$.  In Case I, we choose for each $u\in U_n$
some successor $a(u)$ of $u$ and let $U_{n+1}=\{a(u):u\in U_n\}$.  For
every $a\in U_{n+1}$ we find some $T_a\subseteq(T_n)_a$ and $z_a$ so
that $T_a\Vdash\langle\dx(k):k\in A\cap j_{n+1}\rangle=z_a$, and let
$T_{n+1}=\bigcup\{T_a:a\in U_{n+1}\}$.  In this case $|U_{n+1}|=|U_n|$
and so (vii) holds for $n+1$ as well, while (viii) for $n+1$ follow
either from (viii) or from (vii) for $n$ (the latter if $j_n\in A$).

Thus consider Case II.  For each $u\in U_n$ other than $v_n$ we choose
some $a(u)\in T_n$ of length $l_{n+1}$ such that $a(u)\supset u$,
and find some $T_{a(u)}\subseteq(T_n)_{a(u)}$ and $z_{a(u)}$ so that
$T_{a(u)}\Vdash\langle \dx(k):k\in A\cap m\rangle=z_{a(u)}$.

Let $S$ be the set of all successors of $t$ (which has been chosen so
that $|S|\geq P_m\!^{n+1}$ where $m=\ind(t)$); every $a\in S$ has length
$l_{n+1}$.  For each $a\in S$ we choose $T_a\subseteq(T_n)_a$ and
$z_a$, so that $T_a\Vdash\langle\dx(k):k\in A\cap m\rangle=z_a$.  If
we denote $K=\max(A\cap m)$ then we have
$$
	|\{z_a:a\in S\}|\leq\underset{i\in A\cap m}\to\prod
	N_i\leq\overset{K}\to{\underset{i=0}\to\prod}
	N_i=P_{K+1}\leq P_m,
$$
while $|S|\geq P_m\!^{n+1}$.  Therefore there exists a set $U\subset S$
of size $P_m\!^n$ \st for every $a\in U$ the set $z_a$ is the same.
Therefore if we let
$$
	U_{n+1}=U\cup\{a(u):u\in U_n-\{v_n\}\},
$$
and $T_{n+1}=\bigcup\{T_a:a\in U_{n+1}\}$, $T_{n+1}$ satisfies 
property (iv).  It remains to verify that (vii) and (viii) hold.

To verify (vii), let $k\in A$ be \st $k\geq j_{n+1}=m$.  Since
$m=\ind(t)$, we have $m\notin A$ and so $k>m$.  Let $K\in A$ be
\st $j_n\leq K<m$.  Since $|U_n|<P_K$, we have
$$
	|U_{n+1}|<|U_n|+|U|<P_K+N_m<P_m\cdot N_m=P_{m+1}\leq P_k.
$$

To verify (viii), it suffices to consider only those $k\in A$ \st
$j_n\leq k<m$.  But then $|U_n|<P_k$ and we have
$$
	|\{z_a(k):a\in U_{n+1}\}|\leq|\{z_a:a\in U_{n+1}\}|\leq
	|U_n|+1\leq P_k.
$$  
\line{\hfil \pf}

\enddocument

\end